
\input vanilla.sty
 
\scaletype{\magstep1}
\scalelinespacing{\magstep1}
\def\bull{\vrule height .9ex width .8ex depth -.1ex}

\title Set-functions and factorization
 \endtitle
 
\author N.J.  Kalton and S.J.  Montgomery-Smith\footnote{Both authors
were supported by grants from the National Science Foundation}\\
Department
of Mathematics\\ University of Missouri-Columbia\\ Columbia, Mo. 65211
\endauthor
 
\vskip1truecm
 
\subheading{Abstract} If $\phi$ is a submeasure satisfying an
appropriate lower estimate we give a quantitative result on the total
mass of a measure $\mu$ satisfying $0\le\mu\le\phi.$ We give a dual
result for supermeasures and then use these results to investigate
convexity on non-locally convex quasi-Banach lattices.  We then show how
to use these results to extend some factorization theorems
due to Pisier to the setting of quasi-Banach spaces.  We conclude by
showing
that if $X$ is a quasi-Banach space of cotype two then any operator
$T:C(\Omega)\to X$ is 2-absolutely summing and factors through a Hilbert
space and discussing general factorization theorems for cotype two
spaces.

\vskip2truecm
 
\subheading{1.  Introduction}
 
Let $\Cal A$ be an algebra of subsets of some set $\Omega.$ Let us say
that a set-function $\phi:  \Cal A\to\bold R$ is {\it monotone\/} if it
satisfies $\phi(\emptyset)=0$ and $\phi(A)\le \phi(B)$ whenever
$A\subset B.$ We say $\phi$ is normalized if $\phi(\Omega)=1.$ A
monotone set-function $\phi$ is a {\it submeasure\/} if $$\phi(A\cup
B)\le \phi(A)+\phi(B)$$ whenever $A,B\in\Cal A$ are disjoint, and $\phi$
is a {\it supermeasure\/} if $$\phi(A\cup B)\ge \phi(A)+\phi(B)$$
whenever $A,B\in\Cal A$ are disjoint.  If $\phi$ is both a submeasure
and supermeasure it is a (finitely additive) measure.
 
If $\phi$ and $\psi$ are two monotone set-functions on $\Cal A$ we shall
say that $\phi$ is $\psi$-continuous if $\lim_{n\to\infty}\phi(A_n)=0$
whenever $\lim_{n\to\infty}\psi(A_n)=0.$ If $\phi$ is $\psi-$continuous
and $\psi$ is $\phi-$continuous then $\phi$ and $\psi$ are {\it
equivalent.} A monotone set-function $\phi$ is called {\it exhaustive\/}
if $\lim_{n\to\infty}\phi(A_n)=0$ whenever $(A_n)$ is a disjoint
sequence in $\Cal A.$ The classical (unsolved) Maharam problem ([1],
[5], [6] and [15]) asks whether every exhaustive submeasure is
equivalent
to a measure.  A submeasure $\phi$ is called {\it pathological} if
whenever $\lambda$ is a measure satisfying $0\le\lambda\le \phi$ then
$\phi=0.$ The Maharam problem has a positive answer if and only if there
is no normalized exhaustive pathological submeasure.
 
While the Maharam problem remains unanswered, it is known (see e.g.  [1]
or [15]) that there are non-trivial pathological submeasures.  In the
other direction it is shown in [6] that if $\phi$ is a non-trivial {\it
uniformly exhaustive\/} submeasure then $\phi$ cannot be pathological.
$\phi$ is uniformly exhaustive if given $\epsilon>0$ there exists
$N\in\bold N$ such that whenever $\{A_1,\ldots,A_N\}$ are disjoint sets
in $\Cal A$ then $$\min_{1\le i\le N}\phi(A_i)<\epsilon.$$
 
Let us say that a monotone set-function $\phi$ satisfies an {\it upper
p-estimate} where $0<p<\infty$ if $\phi^p$ is a submeasure, and a {\it
lower p-estimate} if $\phi^p$ is a supermeasure.  If $\phi$ is a
normalized submeasure which satisfies a lower $p$-estimate for some
$1<p<\infty$ then $\phi$ is uniformly exhaustive and hence by results of
[6] there is a non-trivial measure $\lambda$ with $0\le\lambda\le \phi.$
In Section 2 we prove this by a direct argument which yields a
quantitative estimate that $\lambda$ can be chosen so that $$
\lambda(\Omega) \ge 2(2^p-1)^{-1/p}-1.$$ Notice that the expression on
the right tends to one as $p\to 1$ so this result can be regarded for
$p$ close to $1$ as a perturbation result.  The dual result for
supermeasures (Theorem 2.2) is that if a normalized supermeasure $\phi$
satisfies an upper $p$-estimate where $0<p<1$ then there is a
$\phi-$continuous measure $\lambda$ with $\lambda\ge \phi$ and
$$\lambda(\Omega) \le 2(2^p-1)^{-1/p}-1.$$
 
While we believe these results with their relatively simple proofs have
interest in their own right, one of our motivations for considering them
was to use them in the study of some questions concerning quasi-Banach
lattices, or function spaces.
 
It is well-known ([7]) that a Banach lattice $X$ with a (crude) upper
$p$-estimate is $r$-convex for every $0<r<p$ (for the definitions, see
Section 3).  This result does not hold for arbitrary quasi-Banach
lattices [2]; a quasi-Banach lattice need not be $r$-convex for any
$r<\infty.$ However, it is shown in [2] that if $X$ has a crude lower
$q$-estimate for some $q<\infty$ then the result is true.  We provide
first a simple proof of this fact, only depending on the arguments of
Section 2. We then investigate this result further, motivated by the
fact that if $X$ satisfies a strict upper $p$-estimate (i.e. with
constant one) and a strict lower $p$-estimate then $X$ is $p$-convex
(and in fact isometric to an $L_p(\mu)-$space.)  We thus try to estimate
the constant of $r$-convexity $M^{(r)}(X)$ when $0<r<p$ and $X$ has a
strict upper $p$-estimate and a strict lower $q$-estimate where $p,q$
are close.  We find that an estimate of the form $$ \log M^{(r)}(X) \le
c\theta(1+|\log\theta|)$$ where $c=c(r,p)$ and $\theta=q/p-1.$ We show
by example that such an estimate is best possible.  Let us remark that
in the case $r=1<p<q$ the constant $M^{(1)}(X)$ measures the distance
(in the Banach-Mazur sense) of the space $X$ from a Banach lattice.
 
Finally in Section 4 we apply these results to give extensions of some
factorization theorems of Pisier [13] to the non-locally convex setting.
Pisier showed the existence of a constant $B=B(p)$ so that if $X$ is a
Banach space and $T:C(\Omega)\to X$ is bounded satisfying for a suitable
constant $C$ and all disjointly supported functions $f_1,\ldots, f_n\in
C(\Omega)$ $$ \left(\sum_{k=1}^n \|Tf_k\|^p\right)^{1/p} \le C\max_{1\le
k\le n}\|f_k\|$$ then there is a proabability measure $\mu$ on $\Omega$
so that for $f\in C(\Omega)$ $$ \|Tf\| \le BC\|f\|_{L_{p,1}(\mu)}$$
where $L_{p,1}(\mu)$ denotes a the Lorentz space $L_{p,1}$ with respect
to $\mu.$
 
Pisier's approach in [13] uses duality and so cannot be used in
the
case when $X$ is a quasi-Banach space.  Nevertheless the result can be
extended and we prove that if $0<r<1$ there is a constant $B=B(r,p)$ so
that if $X$ is $r$-normable then there exists a probability measure
$\mu$ so that for all $f\in C(\Omega),$ $$ \|Tf\| \le BC
\|f\|_{L_{p,r}(\mu)}.  $$
 
We apply these results to show if $X$ is a quasi-Banach space of cotype
two then any operator $T:C(\Omega)\to X$ is 2-absolutely summing and so
factorizes through a Hilbert space.  We conclude by presenting a dual
result and make a general conjecture that if $X$ and $Y$ are quasi-Banach
spaces such that $X^*$ and $Y$ have cotype two and $T:X\to Y$ is an
approximable linear operator then $T$ factorizes through a Hilbert space.

\vskip.8truecm

\subheading{2.  Submeasures and supermeasures}
 
Let us define for $0<p<\infty,$ $$K_p=\frac2{(2^p-1)^{1/p}}-1.$$ Notice
that for $p$ close to $1$ we have $K_p\sim 1-4(p-1)\log 2$ while for $p$
large we have $K_p\sim p^{-1}2^{-p}.$
 
We now state our main result on submeasures with a lower estimate (see
Section 1 for the definitions).
 
\proclaim{Theorem 2.1}Let $\Cal A$ be an algebra of subsets of $\Omega.$
Suppose that $\phi$ is a normalized submeasure on $\Cal A$, which
satisfies a lower $p$-estimate where $1<p<\infty.$ Then there is measure
$\lambda$ on $\Cal A$ with $0\le\lambda\le \phi$ and $\lambda(\Omega)\ge
K_p.$ \endproclaim
 
\demo{Proof}By an elementary compactness argument we need only prove the
result for the case when $\Omega$ is finite and $\Cal A=2^{\Omega}.$ We
fix such an $\Omega.$
 
Let $\gamma$ be the greatest constant such that, whenever $\phi$ is a
normalized submeasure on $\Omega$ satisfying a lower $p$-estimate then
there is a measure $\lambda$ with $0\le \lambda\le \phi$ and
$\phi(\Omega)\ge \gamma.$ It follows from a simple compactness argument
that there is a normalized submeasure $\phi$ satisfying a lower
$p$-estimate for which this constant is attained; that is if $\lambda$
is a measure with $0\le\lambda\le\phi$ then $\lambda(\Omega)\le \gamma.$
We choose this $\phi$ and then pick an optimal measure $\lambda$ with
$0\le\lambda\le\phi$ and $\lambda(\Omega)=\gamma.$
 
Let $\delta=(2^p-1)^{-1/p}.$ Let $E$ be a maximal subset of $\Omega$ such
that $\lambda(E)\ge \delta \phi(E)$ and let $F=\Omega\setminus E.$
Suppose $A\subset F;$ then $$ \lambda(A)+\lambda(E) =\lambda(A\cup E)
\le \delta\phi(A\cup E) \le\delta(\phi(A) + \phi(E)),$$ and so $$
\lambda(A) \le \delta\phi(A).  \tag 1$$
 
Let $q$ be the conjugate index of $p$, i.e.  $p^{-1}+q^{-1}=1.$ Let
$\nu$ be any measure on $\Cal A$ so that $0\le\nu\le \phi.$ Suppose
$c_1,c_2\ge 0$ are such that $c_1^q+c_2^q=1.$ Consider the measure $$
\mu(A) =c_1\lambda(A\cap E) + c_2\nu(A\cap F).$$ Then for any $A$
$$
\align \mu(A) &\le (c_1^q+c_2^q)^{1/q}(\lambda(A\cap E)^p + \nu(A\cap
F)^p)^{1/p}\\ &\le (\phi(A\cap E)^p +\phi(A\cap F)^p)^{1/p} \le \phi(A).
\endalign $$ Hence $\mu(\Omega)\le \gamma$ which translates as $$
c_1\lambda(E) + c_2\nu(F) \le \gamma.$$ Taking the supremum over all
$c_1,c_2$, we have:  $$ \lambda(E)^p + \nu(F)^p \le \gamma^p.\tag 2$$
 
Now take $\nu(A)=\delta^{-1}\lambda(A\cap F).$ It follows from (1) that
$\nu\le\phi$ and hence from (2), $$ \lambda(E)^p
+\delta^{-p}\lambda(F)^p \le \gamma^p.$$ If we set $t=\lambda(E)/\gamma$
then $$ t^p + (2^p-1)(1-t)^p \le 1$$ and it follows by calculus that
$t\le \frac12.$ Hence $\lambda(E)\le\gamma/2.$
 
Now, consider the submeasure $\psi(A)=\phi(A\cap F).$ By hypothesis on
$\gamma$ there exists a positive measure $\nu$ on $\Omega$ such that
$0\le\nu\le\psi$ and $\nu(\Omega)\ge \gamma \psi(\Omega).$ Thus for all
$A$ we have $0\le \nu(A)\le \phi(A\cap F),$ $\nu(E)=0$ and
$\nu(F)=\nu(\Omega)\ge\gamma \phi(F).$ Returning to equation (2) we
have:  $$ \lambda(E)^p+\gamma^p\phi(F)^p \le \gamma^p.$$ However
$\phi(F)\ge 1-\phi(E)\ge 1-\delta^{-1}\lambda(E).$ Thus, recalling that
$t=\lambda(E)/\gamma$ $$ t^p + (1-\delta^{-1}\gamma t)^p\le 1,$$ which
simplifies to $$ \gamma \ge
\delta\left(\frac{1-(1-t^p)^{1/p}}{t}\right).$$ Since $t\ge \frac12$ it
follows again by calculus arguments that the right-hand side is
minimized when $t=\frac12$ and then $$ \gamma \ge
2\delta\left(1-(1-2^{-p})^{1/p}\right)=K_p$$ and this completes the
proof.   \bull\enddemo
 
In almost the same manner, we can prove the dual statement for
supermeasures.
 
\proclaim{Theorem 2.2}Let $\Cal A$ be an algebra of subsets of a set
$\Omega$.  Suppose $0<p<1$ and that $\phi$ is a normalized supermeasure
on $\Omega$ which satisfies an upper $p$-estimate.  Then there is a
$\phi$-continuous measure $\lambda$ on $\Cal A$ such that
$\lambda\ge\phi$, $\lambda(\Omega)\le K_p$.  \endproclaim
 
\demo{Proof}We first prove the existence of some measure $\lambda$ with
$\lambda\ge\phi$ and $\lambda(\Omega)\le K_p$ without requiring
continuity.  As in the preceding proof it will suffice to consider the
case when $\Omega$ is finite and $\Cal A=2^{\Omega}.$ In this case there
is a least constant $\gamma<\infty$ with the property that if $\phi$ is
a normalized supermeasure on $\Omega$ then there is a measure
$\lambda\ge \phi$ with $\lambda(\Omega)\le\gamma.$ We again may choose
an extremal $\phi$ and associated extremal $\lambda$ for which
$\lambda(\Omega)=\gamma.$
 
Define $\delta=(2^p-1)^{-1/p}>1$ we now let $E$ be a maximal subset so
that $\lambda(E)\le \delta\phi(E)$ and defining $F=\Omega\setminus E$ we
obtain in this case that if $A\subset F$ then $\lambda(A)\ge \phi(A).$
 
In this case let $q$ be defined by $\frac1q=\frac1p-1.$ Let $\nu$ be any
measure on $\Cal A$ such that $\nu(A)\ge \phi(A)$ whenever $A\subset F.$
Suppose $c_1,c_2>0$ satisfy $c_1^{-q}+c_2^{-q}=1.$ Consider the measure
$$ \mu(A) = c_1\lambda(A\cap E) +c_2\nu(A\cap F).$$ Then for any $A$, $$
\align \phi(A) &\le (\phi(A\cap E)^p+\phi(A\cap F)^p)^{1/p}\\ &\le
(\lambda(A\cap E)^p + \nu(A\cap F)^p)^{1/p}\\ &\le (c_1\lambda (A\cap E)
+c_2\nu(A\cap F))(c_1^{-q}+c_2^{-q})\\ &= \mu(A).  \endalign $$ Hence
$\mu(\Omega)\ge \gamma$ and so $$ c_1\lambda(E) + c_2\nu(F) \ge
\gamma.$$ Minimizing over $c_1,c_2$ yields $$ \lambda(E)^p +\nu(F)^p \ge
\gamma^p.\tag 3$$ In particular if we let $\nu(A) =
\delta^{-1}\lambda(A\cap F)$ and set $t=\lambda(E)/\gamma$ we obtain $$
t^p + (2^p-1)(1-t)^p \ge 1.$$ Since in this case $p<1$ we are led to the
conclusion that $t\ge \frac12.$
 
Next we consider the supermeasure $\psi(A)=\phi(A\cap F)$ and deduce the
existence of a measure $\nu\ge \psi$ with $\nu(\Omega)\le \gamma
\psi(\Omega)=\gamma\phi(F).$ In this case (3) gives that $$ \lambda(E)^p
+\gamma^p\phi(F)^p \ge \gamma^p.$$ Now as $\phi(F) \le 1-\phi(E) \le
1-\delta^{-1}\lambda(E)$ we have:  $$ t^p + (1-\delta^{-1}\gamma t)^p
\ge 1$$ and this again leads by simple calculus to the fact that
$\gamma\le K_p.$ This then completes the proof if we do not require
continuity of $\lambda.$
 
Now suppose $\phi$ is a normalized supermeasure on $\Cal A$ satisfying
an upper $p$-estimate.  Let $\lambda$ be a minimal measure subject to
the conditions $\lambda\ge \phi$ and $\lambda(\Omega)\le K_p.$ (It
follows from an argument based on Zorn's Lemma that such a minimal
measure exists).  Suppose
$\lim_{n\to\infty}\phi(F_n)=0.$ Consider the measures
$\lambda_n(A)=\lambda(A\cap E_n)$ where $E_n=\Omega\setminus F_n.$ Let
$\Cal U$ be any free ultrafilter on the natural numbers and define
$\lambda_{\Cal U}(A)=\lim_{\Cal U}\lambda_n(A).$ Clearly $\lambda_{\Cal
U}\le \lambda$.  Now for any $A$ $$ \lambda_n(A) =\lambda(A\cap E_n)\ge
\phi(A\cap E_n)\ge (\phi(A)^p-\phi(F_n))^{1/p}.$$ Hence $\lambda_{\Cal
U}=\lambda$ by minimality.  Thus $\lim_{\Cal U}\lambda(F_n)=0$ for every
such ultrafilter and this means
$\lim_{n\to\infty}\lambda(F_n)=0.$\bull\enddemo
 
The following corollary is proved for more general uniformly exhaustive
submeasures in [6].
 
\proclaim{Corollary 2.3}Let $\Cal A$ be an algebra of subsets of
$\Omega$ and let $\phi$ be a submeasure on $\Omega$ such that for some
constant $c>0$ and some $q<\infty,$ we have:  $$ \phi(A_1\cup\cdots\cup
A_n) \ge c(\phi^q(A_1)+\cdots +\phi^q(A_n))^{1/q}$$ whenever
$A_1,\ldots,A_n$ are disjoint.  Then there is a measure $\mu$ on $\Cal
A$ such that $\mu$ and $\phi$ are equivalent.\endproclaim
 
\demo{Proof}Define $\psi$ by $$\psi(A) =\sup(\sum_{k=1}^n\phi^q(A_k))$$
where the supremum is computed over all $n$ and all disjoint
$(A_1,\ldots,A_n)$ so that $A=\cup_{k=1}^nA_k.$ It is not difficult to
show that $\psi$ is a supermeasure satisfying $1/q-$upper estimate and
clearly $c^q\psi\le\phi^q\le\psi.$ By Theorem 2.2 we can pick a measure
$\mu\ge \psi$ which is equivalent to $\psi$ and hence to
$\phi.$\bull\enddemo \vskip.8truecm

\subheading{3.  Convexity in lattices}
 
Let $\Omega$ be a compact Hausdorff space and suppose $\Cal B(\Omega)$
denotes the $\sigma-$algebra of Borel subsets of $\Omega$.  Let
$B(\Omega)$ denote the space of all real-valued Borel functions on
$\Omega$.  An admissible extended-value quasinorm on $B(\Omega)$ is a
map $f\to \|f\|_X,$ $(B(\Omega)\to [0,\infty])$ such that:\newline (a)
$\|f\|_X\le\|g\|_X$ for all $f,g\in B(\Omega)$ with $|f|\le|g|$
pointwise.\newline (b) $\|\alpha f\|_X = |\alpha|\|f\|_X$ for $f\in
B(\Omega),\ \alpha\in\bold R$\newline (c) There is a constant $C$ so
that if $f,g\ge 0$ have disjoint supports then $\|f+g\|_X\le
C(\|f\|_X+\|g\|_X).$ \newline (d) There exists a strictly positive $u$
with $0<\|u\|_X<\infty.$ \newline (e) If $f_n\ge 0$ and $f_n\uparrow f$
pointwise, then $\|f_n\|_X\to\|f\|_X.$ \newline The space
$X=\{f:\|f\|_X<\infty\}$ is then a quasi-Banach function space on
$\Omega$ equipped with the quasi-norm $\|f\|_X$ (more precisely one
identifies functions $f,g$ such that $\|f-g\|_X=0$).  We say that $X$ is
order-continuous if, in addition, we have:  \newline (f) If
$f_n\downarrow 0$ pointwise and $\|f_1\|_X<\infty$ then
$\|f_n\|_X\downarrow 0.$
 
Conversely if $X$ is a quasi-Banach lattice which contains no copy of
$c_0$ and has a weak order-unit then standard representation theorems
can be applied to represent $X$ as an order-continuous quasi-Banach
function space on some compact Hausdorff space $\Omega$ in the above
sense.  More precisely, if $u$ is a weak order-unit then there is a
compact Hausdorff space $\Omega$ and a lattice embedding $L:C(\Omega)\to
X$ so that $L[0,\chi_{\Omega}]=[0,u]$.  Since $X$ contains no copy of
$c_0$ we can use a result of Thomas [16] to represent $L$ is the form $$
Lf =\int_{\Omega}fd\Phi$$ where $\Phi$ is regular $X$-valued Borel
measure on $\Omega.$ This formula then extends $L$ to all bounded Borel
functions.  We now define the quasi-Banach function space $Y$ by $$
\|f\|_Y = \sup_n\| L(\min(|f|,n\chi_{\Omega})\|_X $$ and it may be
verified by standard techniques that $L$ extends to a lattice
isomorphism of $Y$ onto $X$ (which is an isometry if we assume that the
quasi-norm on $X$ is continuous).
 
For an arbitrary quasi-Banach function space $X$ and $0<p<\infty$ we
define the $p$-convexity constant $M^{(p)}(X)$ to be the least constant
(possibly infinite) such that for $f_1,\ldots,f_n\in X$ $$
\|(\sum_{i=1}^n|f_i|^p)^{1/p}\|_X \le
M^{(p)}(\sum_{i=1}^n\|f_i\|_X^p)^{1/p}$$ and we let the $p$-concavity
constant $M_{(p)}(X)$ be the least constant such that $$
M_{(p)}\|(\sum_{i=1}^n|f_i|^p)^{1/p}\|_X \ge
(\sum_{i=1}^n\|f_i\|_X^p)^{1/p}.$$ We also let $M^{(0)}(X)$ be the least
constant such that $$ \||f_1\ldots f_n|^{1/n}\|_X \le
M^{(0)}(\prod_{i=1}^n\|f_i\|_X)^{1/n}.$$
 
$X$ is called $p$-convex if $M^{(p)}(X)<\infty$ and $p$-concave if
$M_p(X)<\infty;$ we will say that $X$ is geometrically convex if
$M^{(0)}(X)<\infty.$ In [2] $X$ is called L-convex if it is $p$-convex
for some $p>0$; it is follows from [2] and [4] that $X$ is L-convex if
and only if it is geometrically convex.
 
Let us now turn to upper and lower estimates.  We say $X$ satisfies a
crude upper $p$-estimate with constant $a$ if for any disjoint
$f_1,\ldots,f_n$ we have $$ \|f_1+\ldots f_n\|_X \le
a(\sum_{i=1}^n\|f_i\|_X^p)^{1/p}$$ and we say that $X$ satisfies an
upper $p$-estimate if $a=1.$ We say that $X$ satisfies a crude lower
$q$-estimate with constant $b$ if for any disjoint $f_1,\ldots,f_n$ we
have $$ b\|f_1+\ldots f_n\|_X \ge (\sum_{i=1}^n\|f_i\|_X^p)^{1/p};$$ and
$X$ satisfies a lower $q$-estimate if $b=1.$

\proclaim{Lemma 3.1}Suppose $0<p<q<\infty.$ If $X$ is a quasi-Banach
function space satisfying a crude upper $p$-estimate with constant $a$
and a crude lower $q$-estimate with constant $b$ then there is an
equivalent function space quasinorm $\|\,\|_Y$ satisfying an upper $p$
and a lower $q-$estimate with $$ \|f\|_X \le \|f\|_Y \le ab\|f\|_X.$$
\endproclaim

\demo{proof}First we define $$ \|f\|_W =\inf
(\sum_{i=1}^n\|f\chi_{A_k}\|_X^p)^{1/p}$$ where the infimum is taken
over all possible Borel partitions $\{A_1,\ldots,A_n\}$ of $\Omega$.  It
is clear that $\|f\|_W\le \|f\|_X\le a\|f\|_W$ and it can be verified
that $W$ satisfies an upper $p$-estimate and a crude lower $q$-estimate
with constant $b$.  Next we define $$ \|f\|_V =\sup
(\sum_{i=1}^n\|f\chi_{A_k}\|_W^q)^{1/q}$$ and finally set
$\|f\|_Y=a\|f\|_V$.  We omit the details.\bull\enddemo
 
We now give a simple proof of the result proved in [2] that any
quasi-Banach function space which satisfies a lower estimate is
L-convex.  We recall first that if $\mu$ is any Borel measure on
$\Omega$ then $L_{p,\infty}(\mu)$ is the space of all Borel functions
such that $$ \|f\|_{L_{p,\infty}(\mu)} = \sup_{t>0}
t(\mu\{|f|>t\})^{1/p}<\infty.$$
 
\proclaim{Theorem 3.2}Let $X$ be a p-normable quasi-Banach function
space which satisfies a crude lower $q$-estimate.  Then:  \newline (i)
$X$ is $r$-convex for $0<r<p$\newline (ii) There is a measure $\mu$ on
$\Omega$ such that $\|f\|_X=0$ if and only if $f=0$ $\mu-$a.e.
\endproclaim
 
\demo{proof}We may assume by Lemma 3.1 that $X$ has an upper
$p$-estimate and a lower $q$-estimate.  Now suppose $f_1,\ldots,f_n\in
X_+$ and $(\sum_{i=1}^nf_i^r)^{1/r}=f.$ Consider the submeasure
$\phi(A)=\|f\chi_A\|_X^p$ for $A\in \Cal B(\Omega).$ This has a lower
$q/p-$estimate and hence there is a Borel measure $\mu$ with
$\mu(\Omega)\ge K_{q/p}\phi(\Omega)$ and such that $\mu(A) \le
\|f\chi_A\|^p$ for any Borel set $A$.  Thus for any $g\in B(\Omega)$
(note that $\mu$ is supported on the set where $f$ is finite) $$
\|gf^{-1}\|_{L_{p,\infty}(\mu)} \le \|g\|_X.$$ Now the space
$L_{p,\infty}(\mu)$ is $r$-convex with $M^{(r)}(L_{p,\infty}(\mu))\le
C=C(p,r)$.  Thus $$ \align \|f\|_X &= \phi(\Omega)^{1/p} \le
K_{q/p}^{-1/p}\|\chi_{\Omega}\|_{L_{p,\infty}(\mu)} \\ &\le
CK_{q/p}^{-1/p}(\sum_{i=1}^n\|f_if^{-1}\|_{L_{p,\infty}(\mu)}^r)^{1/r}\\
&\le C'(\sum_{i=1}^n\|f_i\|_X^r)^{1/r} \endalign $$ where
$C'=C'(p,q,r).$ For (ii) let $u$ be a strictly positive function with
$0<\|u\|_X<\infty$ and define $\phi(A)=\|u\chi_A\|^p;$ then by Corollary
2.3 there is a measure $\mu$ equivalent to $\phi$ and the conclusion
follows quickly.\bull\enddemo
 
Now suppose $\mu$ is any (finite) Borel measure on $\Omega$.  We define
the Lorentz space $L_{p,q}$ for $0<p,q<\infty$ by $$ \|f\|_{L_{p,q}} =
\left(\int_0^{\infty}
\frac{q}{p}t^{q/p-1}f^*(t)^q dt\right)^{1/q}.$$ Here $f^*$ is the
decreasing rearrangement of $|f|$ i.e. $f^*(t)=\inf_{\mu(E)\le
t}\sup_{\omega\notin E}|f(\omega)|.$  It can easily seen by integration
by parts that $$\|f\|_{L_{p,q}}
=\left(\int_0^{\infty}qt^{q-1}\mu(|f|>t)^{q/p}dt\right)^{1/q}$$.   It is
then clear that if $p\le q$ then $L_{p,q}$ satisfies an upper $p$ and a
lower
$q$-estimate.  If $p>q$ then $L_{p,q}$ has an upper $q$ and a lower
$p-$estimate.
 
Suppose $p<q$.  We define the functional $$ \|f\|_{\Lambda_{p,q}} = \sup
(\sum_{i=1}^n (\inf_{\omega\in
A_i}|f(\omega)|)^{q}\mu(A_i)^{q/p})^{1/q}\tag 4$$ where the supremum is
taken over all Borel partitions $\{A_1,\ldots,A_n\}$ of $\Omega$.
 
\proclaim{Proposition 3.3}Suppose $0<p<q$.  Then:\newline (i) The
$\Lambda(p,q)-$quasi-norm is the smallest admissible quasi-norm which
satisfies a lower $q$-estimate and such that $\|\chi_A\|
\ge\mu(A)^{1/p}$ for any Borel set.  \newline (ii) If $f\in B(\Omega)$
and $f^*$ is the decreasing rearrangement of $|f|$ on $[0,\infty)$ then
$$ \|f\|_{\Lambda_{p,q}} = \sup_{\Cal T}(\sum_{j=1}^n
f^*(\tau_j)^q(\tau_j-\tau_{j-1})^{q/p})^{1/q} \tag 5$$ where $\Cal
T=\{\tau_0=0<\tau_1<\cdots< \tau_n\}$ runs through all possible finite
subsets of $\bold R.$\newline (iii) If $f\in B(\Omega)$ then $$
\|f\|_{\Lambda_{p,q}} \le \|f\|_{L_{p,q}} \le
((1+\theta)^{2(1+\theta)}\theta^{-\theta})^{1/q}\|f\|_{\Lambda_{p,q}}$$
where $1+\theta=\frac{q}{p}.$ \endproclaim
 
\demo{Proof}(i) is clear from the definition.
 
(ii):  Suppose $f\in B(\Omega)$ and let $\{A_1,\ldots,A_n\}$ be any
Borel partition of $\Omega$.  Suppose that $1\le j,k\le n$.  Suppose
$\inf_{A_j}|f| \le \inf_{A_k}|f|$.  Then it is easy to verify that if we
let $A'_k=\{\omega\in A_j\cup A_k:|f(\omega)|\ge \inf_{A_k}|f|\}$ and
let $A'_j=(A_j\cup A_k)\setminus A'_k$ then the partition obtained by
replacing $A_j,A_k$ by $A'_j,A'_k$ increases the right-hand side of (4).
In particular it follows that (5) defines the $\Lambda_{p,q}$ quasinorm
when $\mu$ is nonatomic.  Further if $f^*$ is constant on an interval
$[\alpha,\beta]$ it suffices to consider $\Cal T$ where no $\tau_j$ lies
in $(\alpha,\beta)$ and this yields the conclusion for general $\mu.$

(iii):  The first inequality in (iii) is immediate from (i) since
$L_{p,q}$ satisfies a lower $q$-estimate. For
the right-hand inequality we observe that if $h=1+\theta=q/p$:  $$ \align
\|f\|_{L_{p,q}} &= (\int_0^{\infty} \frac{q}{p}t^{q/p-1}(f^*(t))^{q}dt
)^{1/q}\\ &\le
(\sum_{n=-\infty}^{\infty}\frac{q}{p}(h^{n+1}-h^n)h^{(n+1)(q/p-1)}
f^*(h^n)^{q})^{1/q}\\ &=
(\sum_{n=-\infty}^{\infty}h(h-1)h^{q/p-1}h^{q(n+1)/p}f^*(h^{
n+1})^{q})^{1/q}\\ &\le
h^{\frac2p}(h-1)^{(\frac1q-\frac1p)}\|f\|_{\Lambda_{p,q}}.  \endalign $$
The result then follows.\bull\enddemo

Under the hypothesis $p>q$ we define $\Lambda_{p,q}$ by $$
\|f\|_{\Lambda_{p,q}} = \inf(\sum_{i=1}^n (\sup_{\omega\in
A_i}|f(\omega)|^q)\mu(A_i)^{q/p})^{1/q}$$ where the infimum is again
computed over all Borel partitions of $\Omega.$ Proposition 3.4 now has
an analogue whose proof is very similar and we omit most of the details.

\proclaim{Proposition 3.4}Suppose $0<q<p.$ Then:  \newline (i) The
$\Lambda(p,q)$-quasinorm is the largest admissible quasi-norm which
satisfies an upper $q$-estimate and such that
$\|\chi_A\|\le\mu(A)^{1/p},$ for any Borel set $A.$ \newline (ii) If
$f\in B(\Omega)$ then $$ \|f\|_{\Lambda_{p,q}} = \inf_{\Cal
T}(\sum_{j=1}^n f^*(\tau_{j-1})^q(\tau_j-\tau_{j-1})^{q/p})^{1/q} $$ where
$\Cal T=\{\tau_0=0<\tau_1<\cdots <\tau_n=\mu(\Omega)\}$ runs through all
possible finite subsets of $[0,\mu(\Omega)].$ \newline (iii) If $f\in
B(\Omega)$ then $$ \|f\|_{L_{p,q}} \ge \|f\|_{\Lambda_{p,q}} \ge
((1+\theta)^{-2-\theta}\theta^{\theta})^{1/p}\|f\|_{L_{p,q}} $$ where
$\theta=\frac{p}{q}-1.$ \endproclaim

\demo{Proof} We will only proof the second inequality in (iii).  We
define $h=1+\theta=p/q$.
 
$$ \align \|f\|_{L_{p,q}} &=
(\int_0^{\infty}\frac{q}{p}t^{q/p-1}f^*(t)^qdt)^{1/q}\\ &\ge
(\sum_{n=-\infty}^{\infty}\frac{q}{p}(h^{n+1}-h^n)h^{n(q/p-1)}f^*(h^{(n+1)})^q
)^{1/q}\\ &=
(\sum_{n=-\infty}^{\infty}h^{-1}(h-1)h^{(n-1)q/p}f^*(h^n)^q)^{1/q}\\
&\ge (h-1)^{(\frac1q-\frac1p)}h^{-1/p-1/q}\|f\|_{\Lambda_{p,q}}.
\endalign $$ The result then follows.\bull\enddemo

We now immediately deduce the following:
 
\proclaim{Proposition 3.5}Let $X$ be quasi-Banach function space on
$\Omega$ satisfying a crude upper $p$-estimate with constant $a$ and a
crude lower $q$-estimate with constant $b$.  Then if $f\in X_+$ with
$\|f\|_X=1$:  \newline (i) There is a probability measure $\mu$ on
$\Omega$ such that $f>0$ $\mu-$a.e. and if $g\in X$, $$
\|gf^{-1}\|_{\Lambda_{p,q}(\mu)} \le abK_{q/p}^{-1/p}\|g\|_X.$$ \newline
(ii) There is a probability measure $\lambda$ on $\Omega$ such that
$f>0$ $\lambda-$a.e. and if $g\in X$ $$ \|g\|_X \le
abK_{p/q}^{1/q}\|gf^{-1}\|_{\Lambda_{q,p}}(\lambda).$$ \endproclaim
 
\demo{proof}We first introduce an equivalent quasinorm $\|\,\|_Y$ with
an exact upper $p$ and lower $q$-estimate as in Lemma 3.1 so that
$\|g\|_X\le\|g\|_Y\le ab\|g\|_X$ for all $g$.
 
(i)As in Theorem 3.2 we consider the submeasure
$\phi(A)=\|f\chi_A\|_Y^p.$ There is a probability measure $\mu$ such
that $$0\le\mu\le K_{q/p}^{-1} \frac{\phi(A)}{\phi(\Omega)}$$ for all
Borel sets $A.$ Then for $g\in X,$ and any Borel partition
$\{A_1,\ldots,A_n\}$ of $\Omega,$ $$ \align
(\sum_{i=1}^n(\inf_{A_i}|gf^{-1}|^q)\mu(A_i)^{q/p})^{1/q} &\le
K_{q/p}^{-1/p}\phi(\Omega)^{-1/p} (\sum_{i=1}^n
(\inf_{A_i}|gf^{-1}|^q)\|f\chi_A\|_Y^{q})^{1/q} \\ &\le
K_{q/p}^{-1/p}\|f\|_Y^{-1}\|g\|_Y\\ &\le K_{q/p}^{-1/p}ab\|g\|_X.
\endalign $$ Thus (i) follows.  The proof of (ii) is very similar.  In
this case we consider the supermeasure $\phi(A)=\|f\chi_A\|_Y^q$.  There
is a probability measure $\lambda$ on $\Omega$ such that $$0\le
K_{p/q}^{-1}\frac{\phi(A)}{\phi(\Omega)} \le \lambda$$ for all Borel
sets $A.$ Thus for $g\in X,$ and any Borel partition
$\{A_1,\ldots,A_n\}$ of $\Omega,$ $$ \align \|g\|_X &\le \|g\|_Y\\ &\le
(\sum_{i=1}^n(\sup_{A_i} |gf^{-1}|^p)\|f\chi_{A_i}\|_Y^p)^{1/p}\\ &\le
(\sum_{i=1}^n(\sup_{A_i} |gf^{-1}|^p) \phi(A_i)^{p/q})^{1/p}\\ &\le
K_{p/q}^{1/q}\|f\|_Y (\sum_{i=1}^n (\sup_{A_i}|gf^{-1}|^p)
\lambda(A_i)^{p/q})^{1/p} \endalign $$ and the result
follows.\bull\enddemo

\proclaim{Theorem 3.6}Suppose $0<r<p<\infty$.  Then there is a constant
$c=c(r,p)$ such that if $X$ is a quasi-Banach function space satisfying
an upper $p$-estimate and a lower $q-$estimate where $q/p=1+\theta<2$
then $\log M^{(r)}(X) \le \theta(c+\frac1p|\log\theta|).$ \endproclaim
 
\demo{Proof}We use Proposition 3.5.  Suppose $f_1,\ldots,f_n$ are
nonnegative functions in $X$ with $\|f\|_X=1$ where
$f=(\sum_{i=1}^nf_i^r)^{1/r}.$ Then there is a probability measure $\mu$
on $\Omega$ with $f>0$ a.e. and such that if $g\in X$ $$
\|gf^{-1}\|_{\Lambda_{p,q}} \le K_{q/p}^{-1/p}\|g\|_X.$$ Notice that
$K_{q/p}^{-1}\le e^{c_1\theta}$ for some $c_1$.  We also have $q/p\le
e^{\theta}$ and $(q-r)/(p-r)\le c_2\theta$ where $c_2$ depends only on
$p,r.$ Let $w_i=f_if^{-1}$.
 
We note first that for any $w\ge 0$ in $B(\Omega)$ we have $$ \align
\int w^rd\mu &= \int_0^1 w^*(t)^rdt\\ &\le
(\frac{q}{p})^{-r/q}\|w\|_{L_{p,q}}^r(\int_0^1
t^{-\frac{r(q-p)}{p(q-r)}}dt)^{1-r/q}\\ &\le
\frac{p}{q}(\frac{q-r}{p-r})^{1-r/q}\|w\|_{L_{p,q}}^r\\ &\le
e^{c_3\theta}\|w\|_{L_{p,q}}^r \endalign $$ where $c_3=c_3(r,p).$ Thus
$$ \int w^r d\mu \le \theta^{-r\theta/q}e^{c_4\theta}\|w\|_{
\Lambda_{p,q}}^r$$ where $c_4=c_4(r,p).$ Applying this to the $w_i$ and
summing we have $$ 1 \le
\theta^{-r\theta/p}e^{c_5\theta}\sum_{i=1}^n\|f_i\|_X^r,$$ with $c_5$
depending only on $r,p.$ The result now follows.\bull\enddemo
 
\demo{Example} We show that the estimate in the previous theorem is
essentially best possible.  For convenience we consider the case $p=1$
and $q=1+\theta.$ The conclusion of the theorem is that, for $0<r<p,$
$M^{(r)}(X) \le \exp(F_r(\theta))$ where $$ \lim_{\theta\to
0}\frac{F_r(\theta)}{\theta|\log \theta|}=1.$$ Since $M^{(0)}(X)\le
M^{(r)}(X)$ for all $r>0$ a similar conclusion is attained in the case
$r=0.$ We establish a converse by considering only the case $r=0.$
 
For $\theta>0$ we let $X=X_{\theta}=\Lambda_{1,q}[0,1]$ and we let
$\kappa(\theta)=M^{(0)}(X).$ We will set
$\phi(\theta)=\exp(-|\log\theta|^{1/2})$ so that
$\lim_{\theta\to0}\phi(\theta)=0.$ We further set
$\psi(\theta)=(2^{1/q}-1)^{-2}$; then $\psi(\theta)=1+4\theta\log2
+O(\theta^2).$
 
We will consider the function $f=f_{\theta}\in X$ defined by
$f(t)=t^{-1}\chi_{[1-\phi,1]}.$ It follows easily from the definition of
$M^{(0)}(X)$ that $$ \exp(-\frac1\phi\int_{1-\phi}^1\log t\,
dt)\|\chi_{[1-\phi,1]}\|_X \le \kappa \|f\|_X.$$ To obtain this one
derives the integral version of geometric convexity  and applies it to
suitable rotations of $f$.
 Thus if
$\beta(\theta)=1+(1-\phi)\phi^{-1}\log(1-\phi)$ then $$ e^{\beta}\phi
\le \kappa\|f\|_X.$$
 
Turning to the estimation of $\|f\|_X$ we note that $\|f\|_X^q$ is the
supremum of expressions of the form $$ \sum_{j=1}^n
(\tau_j-\tau_{j-1})^q\tau_j^{-q}$$ where
$1-\phi=\tau_0<\tau_1<\cdots<\tau_n=1.$ Now if $\tau_j>\psi\tau_{j-1}$
it can be checked that this expression is increased by interpolating
$(\tau_j\tau_{j-1})^{1/2}$ into the partition.  We there fore may
suppose that we consider only partitions where $\tau_j\le
\psi\tau_{j-1}.$ In this case we estimate:  $$ \align
(\tau_j-\tau_{j-1})^q\tau_j^{-q} &\le
(\tau_j-\tau_{j-1})(\tau_j-\tau_{j-1})^{\theta}\tau_{j-1}^{-1+\theta}\\
&\le (\psi-1)^{\theta} (\tau_j-\tau_{j-1})\tau_{j-1}^{-1} \endalign $$
and after summing we get the estimate $$ \|f\|_X^q \le
(\psi-1)^{\theta}|\log(1-\phi)|.$$ Thus $$ \kappa \ge
e^{\beta}\phi(\psi-1)^{-\theta/q}|\log(1-\phi)|^{-1/q}.$$
 
Now for small $\theta$ we can estimate $|\log(1-\phi)|\le (1+\phi)\phi.$
Thus $$ \kappa^q \ge
e^{q\beta}\phi^{\theta}(\psi-1)^{-\theta}(1+\phi)^{-1}$$ so that $$
\liminf_{\theta\to 0}\frac{\log\kappa}{\theta|\log\theta|} \ge
\liminf_{\theta\to 0}\frac{(\log\phi-\log(\psi-1))}{|\log\theta|}\ge
1.$$ Thus we conclude from this calculation and the theorem that $$
\log\kappa(\theta) = -\theta\log\theta + o(\theta|\log\theta|)$$ as
$\theta\to 0.$ \bull\enddemo
 
\vskip.8truecm

\subheading{4.  The factorization theorems of Pisier}
 
We next show how the results of Proposition 3 quickly give extensions of
some factorization theorems due to Pisier [13].  Our approach is valid
for quasi-Banach spaces since it does not depend on any duality.

\proclaim{Theorem 4.1}Suppose $0<r\le p<q<\infty.$ Suppose $\Omega$ is a
compact Hausdorff space and that $Y$ is an $r$-Banach space.  Suppose
$T:C(\Omega)\to Y$ is a bounded linear operator.  Then the following
conditions on $T$ are equivalent:  \newline (i) There is a constant
$C_1$ so that for any $f_1,\ldots,f_n$ in $C(\Omega)$ with disjoint
support we have $$ (\sum_{i=1}^n\|Tf_i\|^q)^{1/q} \le C_1 \max_{1\le
i\le n}\|f_i\|.$$ \newline (ii) There is a constant $C_2$ so that for
any $f_1,\ldots,f_n$ we have $$ (\sum_{i=1}^n\|Tf_i\|^q)^{1/q} \le C_2
\|(\sum_{i=1}^n|f_i|^p)^{1/p}\|.$$ \newline (iii) There is a constant
$C_3$ and a probability measure $\mu$ on $\Omega$ such that for all
$f\in C(\Omega),$ $$ \|Tf\| \le C_3 \|f\|^{1-p/q}(\int
|f|^pd\mu)^{1/q}$$ \newline (iv) There is a constant $C_4$ and a
probability measure $\mu$ on $\Omega$ such that for all $f\in
C(\Omega),$ $$ \|Tf\| \le C_4 \|f\|_{L_{q,r}(\mu)}.$$ \endproclaim
 
\demo{proof}Some implications are essentially trivial.  Thus
$(iv)\Rightarrow (iii)$ and $(ii)\Rightarrow (i)$; $(iii)\Rightarrow
(ii)$ is easy and we omit it.  It remains to show $(i)\Rightarrow (iv).$
 
To do this we first notice that if $f_n$ is a sequence of disjointly
supported functions in $C(\Omega)$ then $Tf_n\to 0.$ Thus the theorem of
Thomas [16] already cited allows us to find a regular $Y$-valued Borel
measure $\Phi$ on $\Omega$ so that $$Tf=\int fd\Phi$$ and use this
formula to extend $T$ to the bounded Borel functions
$B_{\infty}(\Omega).$ It is easy to verify the condition (i) remains in
effect for disjointly supported bounded Borel functions.
 
Now we introduce a quasi-Banach function space $Z$ by defining $$
\|f\|_Z = \sup \left(\sum_{i=1}^n \|Tg\|^q \right)^{1/q} ,$$ where the
supremum is over disjoint $g_i \in B_{\infty}(\Omega)$\ with $|g_i|\le
|f|$.  It is immediate that $\|\,\|_Z$ satisfies an upper $r$-estimate
and a lower $q$-estimate.  Also, we see that $\chi_\Omega \in Z$.  Thus
by Proposition 3.5 we can find a probability measure $\mu$ on $\Omega$
so that for some $C,C_4$ $$ \|g\|_Z \le C \|g\|_{\Lambda_{q,r}(\mu)}\le
C_4\|g\|_{L_{q,r}(\mu)}.\bull$$ \enddemo
 
We now prove the companion factorization result for operators on
$L_p$-spaces.

\proclaim{Theorem 4.2}Suppose $0<r\le q< s< \infty.$ Let $(\Omega,\mu)$
be a
$\sigma$-finite measure space.  Let $Y$ be an $r$-Banach
space, and let $T:L_s(\mu)\to Y$ be a bounded linear operator.  Then the
following conditions are equivalent:  \newline (i) There is a constant
$C_1$ so that for any disjoint $f_1,\ldots,f_n$ in $L_s(\mu)$ we have:
$$
(\sum_{i=1}^n \|Tf_i\|^q)^{1/q} \le C_1
\|\sum_{i=1}^n|f_i|\|_{L_s(\mu)}.$$
\newline (ii) There is a constant $C_2$ and a probability measure
$\lambda$ on $\Omega$ so that for any $f\in L_s(\mu)$ and any Borel set
$E$
we have $$ \|T(f\chi_E)\| \le C_2 \|f\|_{L_s(\mu)} \lambda(E)^{1/q-1/s}$$
\newline (iii) There is a constant $C_3$ and a probability measure
$\lambda$ on $\Omega$ so that for any $f\in L_s(\mu)$ and $g\in
B(\Omega)$ with
$|g|\le 1,$ $$ \|T(fg)\| \le C_3 \|f\|_{L_s(\mu)}
\|g\|_{L_{t,r}(\lambda)}$$ where $1/t=1/q-1/s.$ \newline (iv) There
exists $w\in L_1(\mu)$ with $w\ge 0$ and $\int w\,d\mu=1$ and a constant
$C_4$ such that for all $f\in L_s(\mu)$ $$ \|Tf\| \le C_4
\|fw^{-1/s}\|_{L_{q,r}(w\mu)} .$$ \endproclaim
 
\demo{proof}We omit the simple proof that $(ii)\Rightarrow (i)$.  Also
$(iii)\Rightarrow (ii)$\ and $(iv)\Rightarrow(ii)$\ are obvious.
 
We first consider $(i)\Rightarrow (iii).$ For this direction we define a
quasi-Banach function space $Z$ by $$ \|g\|_Z = \sup \left( \sum_{i=1}^n
\|T(f_ig)\|^q \right)^{1/q} ,$$ where the supremum is over disjoint
$f_i\in L_s(\mu)$\ with $\|\sum_{i=1}^n |f_i| \|_s \le 1$\ and $f_ig \in
L_s(\mu)$.  It is clear that $Z$ satisfies an upper $r$-estimate.  We
show that $Z$\ satisfies a lower $t$-estimate, where $1/t = 1/q - 1/s$.
Let us suppose that we have disjoint $g_1,\ldots,g_m \in Z$\ such that
$$ \sum_{k=1}^m \|g_k\|_Z^t > 1 ,$$ and let $g=g_1+\ldots g_m$.  Then
there are $f_{ik}$\ such that $$ \sum_{k=1}^m \left(\sum_{i=1}^n
\|T(f_{ik}g_k)\|^q \right)^{t/q} = 1 $$ and $$ \| \sum_{i=1}^n f_{ik}
\|_s \le 1 .$$ We let $$ \alpha_k = \left(\sum_{i=1}^n
\|T(f_{ik}g_k)\|^q \right)^{t/qs} $$ and $f'_{ik} = \alpha_k f_{ik}$.
Then we have that $$ \sum_{k=1}^m \sum_{i=1}^n \|T(f'_{ik} g)\|^q \ge 1,
$$ and $$ \| \sum_{k=1}^m \sum_{i=1}^n f'_{ik} \|_s^s \le \sum_{k=1}^m
\alpha_k^s \le 1 ,$$ that is, $\| g\|_Z \ge 1 $.
 
We also notice that $\chi_{\Omega}\in Z$.  Thus by Propositions 3.4 and
3.5, there is a probability measure $\lambda$ so that for some $C_3$ and
all bounded $g$ $$ \|g\|_Z \le C_3 \|g\|_{L_{t,r}(\lambda)} .$$
 
Now we show that $ (ii) \Rightarrow (iv) .$ We notice that $\lambda$ is
$\mu-$continuous, and so by the Radon-Nikodym theorem, we can find
$w=d\lambda/d\mu$.  Then for any measurable set $A$ we have:  $$
\|T(\chi_A w^{1/s}\| \le C_2 \|\chi_A w^{1/s}\|_{L_{s}(\mu)}
\lambda(A)^{1/q-1/s} = C_2 \lambda(A)^{1/q}.$$ Now define a function
space $W$ by $$ \|f\|_W =\sup_{|g|\le |f|, g\in
B_{\infty}}\|T(gw^{1/s})\|.$$ Then $W$ has an upper $r$-estimate and $$
\|\chi_A\|_W \le C_2 \lambda(A)^{1/q} .$$ Hence by Proposition 3.4 $$
\|f\|_W \le C_4 \|f\|_{L_{q,r}(\lambda)}$$ for some $C_4$.
\bull\enddemo
 
To conclude the paper let us observe that Theorem 4.1 can be used to
extend other factorization results to non-locally convex spaces.
Let us first recall that an operator $T:X\to Y$ where $X$ is a Banach
space and $Y$ is a quasi-Banach space is called 2-absolutely summing if
there is a constant $C$ so that for $x_1,\ldots,x_n\in X$ we have
$$ (\sum_{i=1}^n \|Tx_i\|^2)^{1/2}\le C\max_{\|x^*\|\le
1}(\sum_{i=1}^n|x^*(x_i)|^2)^{1/2}.$$
A quasi-Banach space $X$ is of cotype $p$ if there is a constant $C$ so
that if $x_1,\ldots,x_n\in X$ then
$$ \text{Ave }_{\epsilon_i=\pm1}\|\sum_{i=1}^n\epsilon_ix_i\| \le
C(\sum_{i=1}^n\|x_i\|^p)^{1/p}.$$
 
For the following theorem for Banach spaces see [12], p. 62.
 
\proclaim{Theorem 4.3}Suppose $\Omega$ is a compact Hausdorff space and
$Y$ is a quasi-Banach space with cotype two.  Suppose $T:C(\Omega)\to Y$
is a bounded operator.  Then $T$ is 2-absolutely summing and hence there
is a probability measure $\mu$ on $\Omega$ and a constant $C$ so that
$\|Tf\|\le C\|f\|_{L_2(\mu)}$ for $f\in C(\Omega).$\endproclaim
 
\demo{Proof}We may assume that $X$ is an $r$-Banach space where
$0<r<1.$  We first note that if
$f_1,\ldots,f_n\in X$ then since
$X$ has cotype two,
$$
\align
  (\sum_{i=1}^n\|Tf_i\|^2)^{1/2} &\le C_1\text
{Ave}_{\epsilon_i=\pm1}\|\sum_{i=1}^n\epsilon_iTf_i\| \\
   &\le C_1\|T\| \text
{Ave}_{\epsilon_i=\pm1}\|\sum_{i=1}^n\epsilon_if_i\|\\
   &\le C_1\|T\| \|\sum_{i=1}^n|f_i|\|
\endalign
$$
 where $C_1$ depends on the cotype two constant of $X.$  Applying Theorem
4.1, we see that there is a probability measure $\nu$ on $\Omega$ and a
constant $C_2$ so that $\|Tf\|\le C_2\|f\|_{L_{2,r}(\nu)}$ for $f\in
C(\Omega).$  In particular it follows that $\|Tf\|\le
C_3\|f\|_{L_4(\nu)}.$   From this we conclude that if $f_1,\ldots, f_n\in
C(\Omega),$  using Khintchine's inequality,
$$
\align
(\sum_{i=1}^n\|Tf_i\|^2)^{1/2} &
\le C_1\text
{Ave}_{\epsilon_i=\pm1}\|\sum_{i=1}^n\epsilon_iTf_i\| \\
&\le
C_4\text{Ave}_{\epsilon_i=\pm1}\|\sum_{i=1}^n\epsilon_if_i\|_{L_4(\nu)}\\
 &\le C_5 \|(\sum_{i=1}^n|f_i|^2)^{1/2}\|_{L_4(\nu)}\\
&\le C_5 \|(\sum_{i=1}^n|f_i|^2)^{1/2}\|_{C(\Omega)}
\endalign
$$
so that $T$ is 2-absolutely summing.  It now follows from the
Grothendieck-Pietsch Factorization Theorem (which applies to
non-locally convex $X$) that there is a probability measure $\mu$ on
$\Omega$ satisfying the conclusions of the Theorem.\bull\enddemo

\demo{Remarks}We conclude with some comments on Theorem 4.3.  We remark
first that the theorem, taking $X=L_1,$ gives a circuitous proof of
Grothendieck's inequality, which is equivalent to the assertion that
every operator from $C(\Omega)$ to $L_1$ is 2-absolutely summing.
 
We also note that there are, by now, many known non-locally convex spaces
of cotype two.  The most natural examples are the spaces $L_p$ when
$0<p<1;$  in this case, Theorem 4.3 is known, and is a consequence of
work of Maurey [8] (see [17], p. 271).  However, more recently Pisier
[14] has
shown that the spaces $L_p/H_p$ have cotype two when $p<1$ and it follows
from work of Xu [18] that the Schatten ideals $S_p$ have cotype two when
$p<1,$ and for such examples Theorem 4.3 is apparently new.
 
We conclude by noting a dual result and then make a conjecture based on
these observations.  First let us recall that if $X$ is a quasi-Banach
space then its dual
$X^*$ defined in the usual way is a Banach space; here $X^*$ need not
separate the points of $X$ and may indeed reduce to zero.  We define the
canonical map (not necessarily injective) $j:X\to X^{**}.$  The Banach
envelope of $X$ is the closure $\hat X$ of $j(X).$
 
We shall say that an operator $T:X\to Y$
is {\it strongly approximable} if $T$ is in the smallest subspace $\Cal
A$ of
$\Cal L(X,Y)$ containing the finite-rank operators and closed under the
pointwise convergence of uniformly bounded nets.
 
The following theorem is essentially known.

\proclaim{Theorem 4.4}Suppose $(\Omega,\mu)$ is a
 $\sigma-$finite measure space.  Let $X$ be a
quasi-Banach
space such that $X^*$ has cotype $q<\infty.$  Then, for $0<p<1,$ there is a
constant $C=C(p,X)$ so that if $T:X\to L_p(\mu)$ is a strongly approximable
operator then there exists $w\ge 0$ with
$\int w^{s}d\mu\le C,$ where $1/s=1/p-1$, and such that
$\|w^{-1}Tx\|_{L_1(\mu)}\le C\|T\|\|x\|$ for $x\in X.$\endproclaim
 
\demo{Remark} If $X$ is a Banach space, then this theorem is due to
Mezrag [9] [10]
with no approximability assumptions on $T.$  If $X$ is not locally convex
this result is essentially proved in [3] and we here show how to obtain
the actual statement from the equivalent Theorem 2.2 of [3].
Note also
that for spaces $X$ with trivial dual, the theorem holds vacuously since the
only strongly approximable operators are identically zero.
 
\demo{Proof}It is shown in Theorem 2.2 of [3] that given $\epsilon>0$
there
exists $C_0(\epsilon)$ so that, for any probability measure $\nu$, if
$T:X\to L_p(\nu)$ is strongly approximable then
there exists a Borel subset $E$ of $\Omega$ so that $\nu(E)>1-\epsilon$ and
$\int_E |Tx|\,d\mu \le C_0\|T\|\|x\|$ for $x\in X.$  Let us fix
$\epsilon=1/2$ and let $C_0=C_0(1/2).$  Suppose $x_1,\ldots,x_n\in X$ and
let $f=\sum_{i=1}^n|Tx_i|.$  If $\|f\|_p>0,$ define
$\nu=\|f\|_p^{-p}|f|^p\mu.$   Consider the operator $S:X\to L_p(\nu)$
defined by $Sx= |f|^{-1}Tx$ (set $Sx(\omega)=0$ when $f(\omega)=0.$)
Then $\|S\| \le \|f\|_p^{-1}\|T\|.$  Thus there is a Borel subset $E$ of
$\Omega$ with $\nu(E)>1/2$ and so that $\int_E|Sx|\,d\nu\le
C_0\|f\|_p^{-1}\|T\|\|x\|$ for $x\in X.$
Now
$$
\align
 \frac12 \le  \nu(E) &= \int_E \sum_{i=1}^n|Sx_i|\, d\nu\\
                &\le  C_0 \|f\|_p^{-1}\|T\| \sum_{i=1}^n\|x_i\|.
\endalign
$$
 and so we obtain an inequality
 $$ (\int (\sum_{i=1}^n |Tx_i|)^pd\mu)^{1/p}\le
2C_0\|T\|\sum_{i=1}^n\|x_i\|,$$
 where $C_1$ depends only on $p,X.$  Now by the factorization results of
Maurey [8] (see [17] p. 264) we obtain the desired
conclusion.\bull\enddemo
 
\proclaim{Theorem 4.5}Let $X$ be a quasi-Banach space such that $X^*$ has
cotype two.  Suppose $\Omega$ is a compact Hausdorff space and $\mu$
is a $\sigma$-finite measure on $\Omega.$  Then for $0<p<1$, there is a
constant $C=C(p,X)$ so that if $T:X\to L_p(\mu)$ is a strongly approximable
bounded operator, then there exists $v\ge 0,$ with
  $\int v^{t}d\mu\le C$ where $1/t=1/p-1/2,$ and such that
$\|v^{-1}Tx\|_{L_2(\mu)} \le C\|T\|\|x\|$ for $x\in X.$
\endproclaim
 
\demo{Proof}By Theorem 4.4 we can find $w\ge 0$ with
 $\int w^{s}d\mu\le C_0$ such that $\|w^{-1}Tx\|_{L_1}\le C_0\|T\|\|x\|.$
Now since $L_1$ is a Banach space we have $\|wTx\|_{L_1}\le
C_0\|T\|\|x\|_{\hat X}$ so that $wT=Sj$ where $S:\hat X\to L_1$ is
bounded with
$\|S\|\le C_0\|T\|.$  Since $\hat X^*=X^*$ has cotype two and $L_1$ has
cotype two we can factor $S$ through a Hilbert space [11] and then apply
Maurey's factorization results to obtain $u\ge 0$ with
$\int u^{2}d\mu\le C_1$ and such that $\|u^{-1}Sx\|_{L_2}\le
C_1\|T\|\|x\|$ for
$x\in X.$  Letting $v=uw$ completes the proof.\bull\enddemo

\demo{Remarks}We discuss a question motivated by Theorems 4.3 and 4.5.
An operator
$T:X\to
Y$ is called approximable if given any compact set $K\subset X$ and
any $\epsilon>0$ there exists a finite-rank operator $F:X\to F$ with
$\|Tx-Fx\|<\epsilon$ for $x\in K.$  Pisier has shown
that if $X,Y$ are Banach spaces such that $X^*$ and $Y$ have cotype two
and $T:X\to Y$ is an approximable linear operator then $T$ factorizes
through a Hilbert space (see [11], [12]).  Does the same result hold if
we
assume $X,Y$ are quasi-Banach spaces?  Theorems 4.3 and 4.5 provide some
evidence that this perhaps is true.
 
\vskip.8truecm
 
\subheading{References}
 
\item{1.} J.P.R.  Christensen and W. Herer, On the existence of
pathological submeasures and the construction of exotic topological
groups, Math.  Ann. 213 (1975) 203-210.
 
\item{2.} N.J.  Kalton, Convexity conditions on non-locally convex
lattices, Glasgow Math.  J. 25 (1984) 141-152.
 
\item{3.} N.J. Kalton, Banach envelopes of non-locally convex spaces,
Canad. J. Math. 38 (1986) 65-86.
 
\item{4.} N.J.  Kalton, Plurisubharmonic functions on quasi-Banach
spaces, Studia Math. 84 (1987) 297-324.
 
\item{5.} N.J.  Kalton, The Maharam problem, Expose 18, Seminaire
Initiation a l'Analyse, 1988/9 (G.  Choquet, G. Godefroy, M. Rogalski,
J. Saint Raymond, editors) Universit\'e Paris VI/VII, 1991.
 
\item{6.} N.J.  Kalton and J.W.  Roberts, Uniformly exhaustive
submeasures and nearly additive set functions, Trans.  Amer.  Math.
Soc. 278 (1983) 803-816.
 
\item{7.} J. Lindenstrauss and L. Tzafriri, {\it Classical Banach
spaces, II, Function spaces } Springer Verlag, Berlin 1979.
 
\item{8.} B. Maurey, Th\'eor\`emes de factorisation pour les op\'erateurs
lin\'eaires \`a valeurs dans un espace $L^p$, Ast\'erisque 11, 1974.
 
\item{9.} L. Mezrag, Th\'eor\`emes de factorisation, et de prolongement
pour les op\'erateurs \`a vauleurs dans les espaces $L^p,$ pour $p<1,$
Thesis, Universit\'e Paris VI, 1984.
 
\item{10.} L. Mezrag, Th\'eor\`emes de factorisation, et de prolongement
pour les op\'erateurs \`a vauleurs dans les espaces $L^p,$ pour $p<1,$
C. R. Acad. Sci. Paris 300 (1985) 299-302.

\item{11.} G. Pisier, Une th\'eor\`eme sur les op\'erateurs entre espaces
de Banach qui se factorisent par un espace de Hilbert, Ann. \'Ecole Norm.
Sup. 13 (1980) 23-43.
 
\item{12.} G. Pisier, {\it Factorization of linear operators and geometry
of
Banach spaces,} NSF-CBMS Regional Conference Series No. 60, Providence,
1986.
 
\item{13.} G. Pisier, Factorization of operators through $L_{p,\infty}$
or $L_{p,1}$ and non-commutative generalizations, Math.  Ann. 276 (1986)
105-136.
 
\item{14.} G. Pisier, A simple proof of a theorem of Jean Bourgain, to
appear.

\item{15.} M. Talagrand, A simple example of a pathological submeasure,
Math.  Ann. 252 (1980) 97-102.
 
\item{16.} G.E.F.  Thomas, On Radon maps with values in arbitrary
topological vector spaces and their integral extensions, unpublished
paper, 1972.
 
\item{17.} P. Wojtaszczyk, {\it Banach spaces for analysts}, Cambridge
University Press, 1991.

\item{18.}Q. Xu, Applications du th\'eor\`emes de factorisation pour des
fonctions \`a valeurs operateurs, Studia Math. 95 (1990) 273-292.

      \bye